\documentclass{ifacconf}

\usepackage{amsmath} 
\usepackage{amssymb}
\usepackage{amsfonts}
\usepackage{dsfont}
\usepackage{natbib}        

\newcommand{\bq}{\begin{equation}}
\newcommand{\eq}{\end{equation}}
\newcommand{\mR}{\mathbb{R}}

\newcommand{\mC}{\mathbb{C}}

\newcommand{\pperp}{\perp \!\!\!\perp}

\DeclareMathOperator{\rank}{rank}

\newtheorem{example}[thm]{Example}

\newtheorem{proposition}[thm]{Proposition}
\newtheorem{remark}[thm]{Remark}

\begin{document}
\begin{frontmatter}

\title{Differential operator Dirac structures} 


\author{Arjan van der Schaft,} 
\author[Second]{Bernhard Maschke} 

\address{Bernoulli Institute for Mathematics, Computer Science and AI, \\
Jan C. Willems Center for Systems and Control, \\University of Groningen, the Netherlands \\(e-mail: a.j.van.der.schaft@rug.nl)}
\address[Second]{Universit\'e Claude Bernard Lyon 1, CNRS, LAGEP, France \\(e-mail: bernhard.maschke@univ-lyon1.fr)}

\begin{abstract}                
As shown in earlier work, skew-adjoint linear differential operators, mapping efforts into flows, give rise to Dirac structures on a bounded spatial domain by a proper definition of boundary variables. In the present paper this is extended to pairs of linear differential operators defining a formally skew-adjoint relation between flows and efforts. Furthermore it is shown how the underlying repeated integration by parts operation can be streamlined by the use of two-variable polynomial calculus. Dirac structures defined by formally skew adjoint operators and differential operator effort constraints are treated within the same framework. Finally it is sketched how the approach can be also used for Lagrangian subspaces on bounded domains.
\end{abstract}


\end{frontmatter}

\section{Introduction}
Dirac structures are an essential ingredient of port-Hamiltonian systems theory; see e.g. \cite{AEU}, \cite{vds14}, \cite{passivitybook}. They capture the power-conserving interconnection structure in port-based modeling of complex physical systems. Dirac structures for the modeling of distributed parameter physical systems with bounded spatial domain were first introduced in \cite{JGP}. In the case of two coupled balance laws the relation between the distributed flows and efforts is given by a first-order differential operator. Using Stokes' theorem this leads to the definition of {\it boundary} flows and efforts. The resulting Dirac structure, called the Stokes-Dirac structure, involves the flows and efforts on the spatial domain {\it and} the boundary flows and efforts on the boundary of the spatial domain. In the case of {\it linear} distributed parameter systems and a one-dimensional spatial domain consisting of a finite interval the Stokes-Dirac structure was generalized in \cite{legorrec}. Starting from a general skew-adjoint linear differential operator, which is mapping the distributed efforts to the distributed flows, boundary flows and efforts are defined through repeated integration by parts.

The present paper generalizes the set-up considered in \cite{legorrec} to Dirac structures defined by {\it pairs} of linear differential operators relating the distributed flows and efforts. Furthermore, instead of relying on 'brute force' repeated integration by parts, we utilize two-variable polynomial calculus to define the boundary variables. This calculus was originally introduced in \cite{QDF} for optimization and dissipativity analysis of finite-dimensional linear systems. It was also used for spectral factorization \cite{trentelman}, while in \cite{rapisarda} it was shown how this calculus admits to compute {\it state maps} from higher-order differential equations (in time) in input and output variables. In the present paper the {\it time} variable is replaced by the (scalar) {\it spatial} variable, and {\it state maps} become {\it boundary maps}, defining the boundary variables of the resulting Dirac structure.
Furthermore it is shown how this construction can be extended to the definition of differential operator Lagrangian subspaces on spatial domains with boundary. In \cite{journalpaper} this will be combined with the definition of Dirac structures on spatial domains with boundary for the formulation of an extended class of boundary control port-Hamiltonian systems.

\section{Recall of two-variable polynomial matrix calculus}
\label{sec:pol}
In this section we will recall some relevant basics of two-variable polynomial calculus from \cite{QDF}; see also \cite{trentelman, rapisarda}.
A $p \times q$ two-variable polynomial matrix $\Phi(\zeta, \eta)$ is an expression in two indeterminates $\zeta$ and $\eta$ of the form
\begin{equation}\label{eq:rem}
\Phi(\zeta, \eta) := \sum_{k,l=0}^{M} \Phi_{k,l} \zeta^k \eta^l 
\end{equation}
for certain $p \times q$ matrices $\Phi_{k,l}$.
The infinite matrix $\widetilde{\Phi}$ whose $(k,l)$-th block is the matrix $\Phi_{k,l}$, $k,l=0,\ldots,M,$ and is zero everywhere else, is called the {\it coefficient matrix} of $\Phi(\zeta, \eta)$.
Associated to $\Phi(\zeta, \eta)$ and its coefficient matrix $\widetilde{\Phi}$ is the {\it bilinear differential operator}\footnote{Called a bilinear differential form in \cite{QDF}.} $D_{\Phi}$ defined as
\begin{equation}
D_{\Phi}(v,w)(z) = \sum_{k,l =0}^{M} \left[\frac{d^k}{dz^k} v(z)\right]^\top \Phi_{k,l}\frac{d^l}{dz^l} w(z),
\end{equation}
acting on vector-valued functions $v: \mR \to \mR^p, w: \mR \to \mR^q$. 
A $p \times p$ two-variable polynomial matrix $\Phi(\zeta, \eta)$ is called {\it symmetric} if $\Phi(\zeta, \eta) = \Phi^\top (\eta, \zeta)$, or equivalently its coefficient matrix $\widetilde{\Phi}$ is symmetric. A symmetric two-variable polynomial matrix $\Phi(\zeta, \eta)$ defines the quadratic differential operator $D_{\Phi}(v,v)(z)$, cf. \cite{QDF}.

An important fact in the calculus of bilinear differential operators is the following, cf. \cite{QDF}. The {\it derivative} of a bilinear differential operator $D_{\Phi}$ defines another bilinear differential operator
\begin{equation}
\frac{d}{dz} \left(D_{\Phi}(v,w)\right)(z) =:D_{\Psi}(v,w)(z) \, ,
\end{equation}
where $\Psi(\zeta,\eta)$ is the $p \times q$ two-variable polynomial matrix
\begin{equation}
\label{formula}
\Psi(\zeta,\eta) = (\zeta + \eta)\Phi(\zeta,\eta)
\end{equation}
This equality is the two-variable polynomial matrix version of the product rule of differentiation, or (in integral form) of 'integration by parts'.

Any $p \times q$ two-variable polynomial matrix $\Phi(\zeta, \eta)$ can be {\it factorized} as a product $\Phi(\zeta, \eta)= X(\zeta)^\top Y(\eta)$ of (single-variable) polynomial matrices $X(\zeta)$ (of dimension $k \times p$) and $Y(\eta)$ (of dimension $k \times q$). Such a factorization corresponds in a one-to-one manner to a factorization $\tilde{\Pi} = \widetilde{X}^\top\widetilde{Y}$ of the coefficient matrix $\widetilde{\Phi}$ of $\Phi(\zeta, \eta)$. Here $\widetilde{X}, \widetilde{Y}$ are the coefficient matrices of $X(\zeta)$ and $Y(\eta)$: let $X(\zeta)= \sum_{k=0}^M X_k \zeta^k$ then $\widetilde{X}$ is the infinite row matrix with $k$-th block given by $X_k, k=0,1, \cdots, M$, and similarly for $\widetilde{Y}$.
\begin{proposition}\label{prop:canfact}
Let $\Phi\in\mR^{p\times q}(\zeta,\eta)$, with $\widetilde{\Phi}$ its coefficient matrix. Any factorization $\Phi(\zeta,\eta)=X(\zeta)^\top Y(\eta)$ corresponds to a factorization $\widetilde{\Phi}=\tilde{X}^\top \widetilde{Y}$, where $\widetilde{X}, \widetilde{Y}$ are the coefficient matrices of $X(\zeta)$, respectively $Y(\eta)$.
\end{proposition}
Factorizations which correspond to the minimal value $k=\rank(\widetilde{\Phi})$, are called {\it minimal}. They are unique up to premultiplication by a constant nonsingular matrix (\cite{QDF}, \cite{trentelman}).

\section{Differential operator Dirac structures over the infinite spatial domain}

Consider a finite-dimensional linear space $\mathcal{F}:= \mathbb{R}^m$; called the {\it flow space}, with elements denoted by $f$. Consider the dual space $\mathcal{E}=\mathcal{F}^*$; called the {\it effort space}, with elements denoted by $e$. Denote the duality product between $\mathcal{F}$ and $\mathcal{E}$ by $<e|f>$. Identifying $f\in \mathcal{F}$ and $e\in \mathcal{E}$ with vectors $f,e \in \mR^m$ obviously $<e|f> = e^\top f$. Then the {\it bond space} $\mathcal{B}= \mathcal{F} \times \mathcal{E}$ is endowed with the symmetric bilinear form 
\begin{equation}
\label{plus}
<(f_1,e_1),(f_2,e_2)>:= <e_1|f_2> +<e_2|f_1>\; ,
\end{equation}
which has the matrix representation
\begin{equation}\label{eq:Qe}
Q_e = \begin{bmatrix} 0_m & I_m \\[2mm] I_m & 0_m \end{bmatrix}\; .
\end{equation}
Note that $Q_e$ has singular values $+1$ and $-1$, both with multiplicity $m$. Hence the bilinear form $<\cdot,\cdot>$ is {\it indefinite} and {\it symmetric}, as well as {\it non-degenerate}. 

Furthermore, $<\cdot,\cdot>$ gives rise to the following symmetric bilinear form on the set $\mathcal{C}^{\infty}(\mathbb{R}, \mathcal{F} \times \mathcal{E})$ of smooth functions $(f,e): \mathbb{R} \to \mathcal{F} \times \mathcal{E}$ with compact support:
\begin{equation}\label{+}
\begin{array}{l}
\langle \langle (f_1,e_1), (f_2,e_2)\rangle \rangle := \\[2mm]
 \qquad  \int_{- \infty}^{\infty}  <(f_1(z),e_1(z)),(f_2(z),e_2(z))>dz
 \end{array}
\end{equation}
This is a non-degenerate form on $\mathcal{C}^{\infty}(\mathbb{R}, \mathcal{F} \times \mathcal{E})$, in the sense that if $\langle \langle (f_1,e_1), (f_2,e_2) \rangle \rangle=0$ for all compact support $(f_1,e_1)$, then $(f_2,e_2)=0$. A subspace $\mathcal{D}$ of the space of functions $(f,e): \mathbb{R} \to \mathcal{F} \times \mathcal{E}$ of compact support is called a {\it Dirac structure} if $\mathcal{D}=\mathcal{D}^{\pperp}$ where ${}^{\pperp}$ denotes the orthogonal companion with respect to the form $\langle \langle  \cdot,\cdot \rangle \rangle $. Equivalently, $\mathcal{D}$ is a Dirac structure if $\langle \langle  \cdot,\cdot \rangle \rangle $ is zero restricted to $\mathcal{D}$, and moreover $\mathcal{D}$ is maximal with respect to this property. 

We will show that such infinite dimensional Dirac structures can be generated by pairs of linear {\it differential operators}\footnote{An abstract algebraic theory about Dirac structures defined by pairs of differential operators can be found in \cite{dorfman}.}. Indeed, consider systems of differential equations over the spatial domain $\mR$ given by
\bq
\label{InfDir}
F(\frac{d}{dz}) f(z) + E(\frac{d}{dz}) e(z) = 0, \quad z \in \mR,
\eq
where $F(\frac{d}{dz})$ and $e(\frac{d}{dz})$ are square linear differential operators. Denote the corresponding $m \times m$ polynomial matrices by $F(s)$ and $E(s)$ in the indeterminate $s$. Now let $F(s)$ and $E(s)$ satisfy
\bq
\label{skew}
\begin{bmatrix} F(-s) & E(-s) \end{bmatrix} \begin{bmatrix} 0_m & I_m \\[2mm] I_m & 0_m \end{bmatrix}\begin{bmatrix} F^\top (s) \\[2mm] E^\top (s) \end{bmatrix} =0 , \quad \forall s \in \mC,
\eq
together with the maximal rank assumption
\bq
\label{rank}
\rank \begin{bmatrix} F(-s) & E(-s) \end{bmatrix} = m,  \quad \forall s \in \mC
\eq
It follows from \cite{QDF, rapisarda}\footnote{In these references the spatial variable $z \in \mR$ is replaced by the time variable $t \in \mR$.} that the set $\mathcal{S}_{\mathcal{D}}$ of solutions with compact support of \eqref{InfDir} satisfies
\bq
\int_{- \infty}^{\infty}  <(f_1(z),e_1(z)),(f_2(z),e_2(z))>dz=0
\eq
for all $(f_1,e_1), (f_2,e_2) \in \mathcal{S}_{\mathcal{D}}$. Hence $\mathcal{S}_{\mathcal{D}} \subset \mathcal{S}_{\mathcal{D}}^{\pperp}$. Furthermore, assuming for the moment that the rational matrix $G(s):= - E^{-1}(s)F(s)$ is {\it proper}, it follows from standard linear system realization theory that the behavior $\mathcal{S}_{\mathcal{D}}$ is generated by a minimal 'state' space system of the form
\bq
\label{stateD}
\begin{array}{rcl}
\frac{d}{dz} b(z) & = & Ab(z) + Bf(z) \\[2mm]
e(z) & = & Cb(z) + Df(z) 
\end{array}
\eq
with 'state' vector $b(z) \in \mR^n$, where $n$ is the McMillan degree of $G(s)$. Here 'state' has been denoted with quotation marks, since actually $b$ corresponds to the {\it boundary} vector of the behavior $\mathcal{S}$ as we will see in the next section. (In ordinary realization theory the {\it spatial} variable $z$ is replaced by the {\it time} variable $t$.) Then, following the techniques exploited in \cite{2009} (again in a time-domain setting), it follows that actually $\mathcal{S}_{\mathcal{D}}^{\pperp} \subset \mathcal{S}_{\mathcal{D}}$, thus yielding equality $\mathcal{S}_{\mathcal{D}}^{\pperp} = \mathcal{S}_{\mathcal{D}}$, proving that the solution set $\mathcal{S}_{\mathcal{D}}$ of \eqref{InfDir} is actually a Dirac structure. 
\begin{proposition}
Consider $F(s),E(s)$ satisfying \eqref{skew} and \eqref{rank}. Let $\mathcal{S}_{\mathcal{D}}$ be the set of smooth solutions of \eqref{InfDir} with compact support within the spatial domain $\mR$. Then $\mathcal{S}_{\mathcal{D}}$ is a Dirac structure with respect to the bilinear form \eqref{+}.
\end{proposition}
Note furthermore that by \eqref{skew} the 'transfer' matrix $G(s)$ satisfies 
\bq
G(s)=-G^\top (-s)
\eq
Therefore the system \eqref{state} is cyclo-lossless in the sense that there exists an invertible matrix $\Sigma=\Sigma^\top$ satisfying
\bq
A^\top \Sigma + \Sigma A = 0, \quad B^\top \Sigma = C, \quad D=-D^\top
\eq
Hence by defining the skew-symmetric matrix $J:=A\Sigma^{-1}$ the system \eqref{state} can be rewritten into port-Hamiltonian form
\bq
\label{statePH}
\begin{array}{rcl}
\frac{d}{dz} b(z) & = & J\Sigma \, b(z) + Bf(z) \\[2mm]
e(z) & = & B^\top \Sigma \, b(z) + Df(z) 
\end{array}
\eq
with Hamiltonian $\frac{1}{2}b^\top \Sigma b$.

In case $G(s)$ is {\it not} proper, the same results continue to hold. This follows from the fact, see \cite{crouch, vds14} for the constant case, that there always exists a partitioning $\{1, \cdots,n\}=I \cup I^c$ such that by defining the polynomial matrix $\widehat{F}(s)$ as the matrix with $k$-th column given by the $k$-th column of $F(s)$ whenever $k \in I$, and equal to the $k$-th column of $E(s)$ whenever $k \in I^c$, and similarly $\widehat{E}(s)$ as the matrix with $k$-th column given by the $k$-th column of $E(s)$ whenever $k \in I$, and equal to the $k$-th column of $F(s)$ whenever $k \in I^c$, then $\widehat{G}(s)= - \widehat{E}^{-1}(s)\widehat{F}(s)$ is proper. Hence by realization theory we obtain a system as in \eqref{state} with inputs $f_k,k\in I, e_k, k \in I^c,$ and outputs $e_k,k\in I, f_k, k \in I^c$.

\section{Differential operator Dirac structures over a bounded interval }\label{sec:one-port_pass}
In this section we will show how the Dirac structure $\mathcal{S}_{\mathcal{D}}$ over the infinite spatial domain $\mR$ turns into a new Dirac structure when restricting to any finite spatial domain $[\alpha, \beta] \subset \mR$, involving {\it boundary variables} at the end points $\alpha, \beta$. This idea was first introduced in the context of Stokes-Dirac structures on a bounded domain (for arbitrary dimensions) in \cite{JGP}. Subsequently this was extended in \cite{legorrec} to Dirac structures on a finite interval induced by skew-adjoint linear differential operators. The present paper extends this to Dirac structures induced by {\it pairs} of linear differential operators $F(\frac{d}{dz}), E(\frac{d}{dz})$ on a finite interval, satisfying \eqref{skew}, \eqref{rank}. Furthermore, the treatment is extended and simplified by the use of two-variable polynomial matrix calculus. In fact, the mathematics employed in this section is similar to the one in \cite{rapisarda}, where the spatial variable $z$ is replaced by the time variable $t$.

Consider the differential operators $F(\frac{d}{dz})$ and $E(\frac{d}{dz})$ satisfying \eqref{skew} and \eqref{rank}, defining a differential operator Dirac structure over the real line $\mR$. Obviously \eqref{skew} is the same as 
\bq
F(-s)E^\top (s) + E(-s) F^\top (s) =0
\eq
It follows that the two variable expression $F(\zeta)E^\top (\eta) + E(\zeta) F^\top (\eta)$ is zero for $\zeta + \eta=0$, and thus, cf. \cite{QDF, rapisarda},
\bq
\label{fac}
F(\zeta)E^\top (\eta) + E(\zeta) F^\top (\eta) = (\zeta + \eta) \Pi(\zeta,\eta)
\eq
for some two-variable symmetric polynomial matrix $\Pi(\zeta,\eta)$. By using the theory of factorization of two-variable polynomial matrices, cf. Section \ref{sec:pol}, it follows that we can write 
\bq
\label{fac1}
\Pi(\zeta,\eta)= Z^\top (\zeta) \Sigma Z(\eta)
\eq
for some polynomial matrix $Z(s)$ and invertible symmetric matrix $\Sigma$. Hence, following the developments in \cite{rapisarda}, the differential operator $Z(\frac{d}{dz})$ defines a minimal {\it boundary map} such that
\bq
\label{fac2}
Z(\frac{d}{dz}) \ell (z) = b(z), \quad \begin{bmatrix} f(z) \\[2mm] e(z) \end{bmatrix} = \begin{bmatrix} F^\top (\frac{d}{dz}) \\[2mm]  E^\top (\frac{d}{dz})\end{bmatrix} \ell (z),
\eq
where the vector of boundary variables $b(z)$ satisfies \eqref{statePH} with $\Sigma$ determined by \eqref{fac}. Note that $\ell(z)$ is a vector of latent variables parametrizing the elements of the differential operator Dirac structure defined by $F(\frac{d}{dz})$ and $E (\frac{d}{dz})$.

\begin{proposition}
\label{prop:dirac}
Consider $F(\frac{d}{dz}),E(\frac{d}{dz})$ satisfying \eqref{skew}, \eqref{rank}. Consider the factorization \eqref{fac}, \eqref{fac1}, \eqref{fac2}.\\
Then for any $\alpha, \beta \in \mR$ the space of $f(z),e(z),b(\alpha),b(\beta)$, which are solutions of 
\bq
\begin{array}{l}
F(\frac{d}{dz}) f(z) + E(\frac{d}{dz}) e(z) = 0, \quad z \in [\alpha, \beta], \\[3mm]
\begin{bmatrix} f(\alpha) \\[2mm] e(\alpha) \end{bmatrix} = \begin{bmatrix} F^\top (\frac{d}{dz}) \\[2mm]  E^\top (\frac{d}{dz})\end{bmatrix} \ell (\alpha), \;
\begin{bmatrix} f(\beta) \\[2mm] e(\beta) \end{bmatrix} = \begin{bmatrix} F^\top (\frac{d}{dz}) \\[2mm]  E^\top (\frac{d}{dz})\end{bmatrix} \ell (\beta)\\[6mm]
b(\alpha) = Z(\frac{d}{dz}) \, \ell (\alpha), \; b(\beta)= Z(\frac{d}{dz}) \, \ell (\beta) ,
\end{array}
\eq
defines a Dirac structure with respect to the bilinear form
\bq
\begin{array}{l}
\int_{\alpha}^{\beta}  <(f_1(z),e_1(z)),(f_2(z),e_2(z))>dz  \\[3mm]
\qquad -\,  b^\top_1 (\beta) \Sigma \, b_2(\beta) \, + \,b^\top_1 (\alpha) \Sigma \, b_2(\alpha)
\end{array}
\eq
\end{proposition}
Note that by using two-variable polynomial calculus we do not have to rely on the 'state' space realization of the 'transfer' matrix $G(s)$, but we may directly construct the 'state' $b(z)$ and the matrix $\Sigma$ in \eqref{statePH} from the factorization \eqref{fac}; mirroring the treatment (in the time-domain!) in \cite{rapisarda}. In particular, we do not have to rely on $G(s)$ (or the modified version $\widehat{G}(s)$) being proper.

A special case of a differential operator Dirac structure is provided by formally {\it skew-adjoint} differential operators, corresponding to $E(s)=I_m$ (and thus $e(z)=\ell(z)$), and 
\bq
F^\top (s)=-F(-s),
\eq
\begin{example}[Stokes-Dirac structure]
The simplest example of a differential operator Dirac structure is the one defined by the formally skew-adjoint differential operator
\bq 
F(s) = \begin{bmatrix} 0 & s \\[2mm] s & 0 \end{bmatrix}, \; E(s) = \begin{bmatrix} 1 & 0 \\[2mm] 0 & 1 \end{bmatrix}
\eq
In this case \eqref{fac} amounts to
\bq
F(\zeta) + F^\top (\eta) = (\zeta + \eta) \begin{bmatrix} 0 & 1 \\ 1 & 0 \end{bmatrix},
\eq
with $\Sigma= \begin{bmatrix} 0 & 1 \\1 & 0 \end{bmatrix}$ and boundary map $Z(s)= \begin{bmatrix} 1 & 0 \\0 & 1 \end{bmatrix}$. This is the scalar spatial domain version of the Stokes-Dirac structure introduced in \cite{JGP}.
\end{example}
\begin{remark}
The Stokes-Dirac structure as introduced in \cite{JGP} is defined on spatial domains of arbitrary dimension, replacing the scalar spatial differentiation $\frac{d}{dz}$ by the exterior derivative $d$. Similarly, much of the theory developed in the present paper can be extended to higher-dimensional spatial domains by replacing $\frac{d}{dz}$ by $d$, and considering polynomials in $d$.
\end{remark}

A very important case in the definition of the boundary variables $b(z)$ occurs if $\Sigma$ has {\it as many positive singular values as negative singular values}. In this case we can always take $\Sigma$ to be in the canonical form
\bq
\Sigma = \begin{bmatrix} 0 & I_p \\[2mm] I_p & 0 \end{bmatrix}, \quad 2p=n
\eq
Then by denoting
\bq
\begin{bmatrix} f_{\delta} \\ e_{\delta} \end{bmatrix} := b(z)= Z(\frac{d}{dz}) \ell (z), \; f_{\delta}, e_{\delta} \in \mR^p,
\eq
we have 
\bq
\label{balance}
e_1^\top f_2 + e_2^\top f_1 = \frac{d}{dz} \left(e_{\delta 1}^\top f_{\delta 2} + e_{\delta 2}^\top f_{\delta 1}\right)
\eq
Equation \eqref{balance} has an immediate interpretation in terms of {\it power balance}. Indeed, in integral form it amounts
\bq
\begin{array}{l}
\int_{\alpha}^{\beta} e_1^\top (z)f_2(z) + e_2^\top (z)f_1(z) dz = \\[4mm]
e_{\delta 1}^\top (\beta)f_{\delta 2}(\beta) \!+ \! e_{\delta 2}^\top (\beta)f_{\delta 1}(\beta) \! - \! e_{\delta 1}^\top (\alpha)f_{\delta 2}(\alpha) \! - \! e_{\delta 2}^\top (\alpha)f_{\delta 1}(\alpha)
\end{array}
\eq
for any interval $[\alpha, \beta] \subset \mR$. In particular, by taking $f:=f_1=f_2, e:=e_1=e_2,f_{\delta}:=f_{\delta 1}=f_{\delta 2}, e_{\delta}:=e_{\delta 1}=e_{\delta 2}$, this implies the power balance
\bq
\int_{\alpha}^{\beta} e^\top (z)f(z) dz = e_{\delta }^\top (\beta)f_{\delta }(\beta) - e_{\delta }^\top (\alpha)f_{\delta }(\alpha)
\eq
where the left-hand side is the total incoming power on the interval $[\alpha, \beta]$, and the right-hand side is the difference of the outgoing power at the right-end point $\beta$ and the outgoing power at the left-end point $\alpha$. The variables $f_{\delta }, e_{\delta }$ are called boundary {\it power} variables, since their product $e_{\delta}^\top f_{\delta}$ equals power.

The condition that $\Sigma$ has as many positive as negative singular values can be verified as follows. Consider the two-variable polynomial matrix $\Pi(\zeta,\eta)$ in \eqref{fac} obtained from dividing $F(\zeta)E^\top (\eta) + E(\zeta) F^\top (\eta)$ by $\zeta + \eta$. Then 
\begin{proposition}[\cite{trentelman}](Prop. 2.1)
$\Sigma$ has as many positive as negative singular values if and only if the coefficient matrix $\widetilde{\Pi}$ of $\Pi(\zeta,\eta)$ has.
\end{proposition}

\begin{remark}
This should be compared with the approach taken in \cite{legorrec}, where boundary variables similar to $f_{\delta}, e_{\delta}$ are defined {\it in general}. However, this is enforced at the expense of 'mixing' the values at the left boundary $\alpha$ and the right boundary $\beta$; in this way extending the matrix $\Sigma$ to a matrix of double dimension
\bq
\begin{bmatrix} \Sigma & 0 \\ 0 & - \Sigma \end{bmatrix},
\eq
which obviously has as many positive and negative singular values.
\end{remark}

The boundary map $b(z)=Z(\frac{d}{dz}) \ell(z)$ can be given the following interpretation. Consider $f_1(z),e_1(z)$ satisfying \eqref{InfDir} on an interval $[\alpha, \gamma]$ and $f_2(z),e_2(z)$ satisfying \eqref{InfDir} on an interval $[\gamma, \beta]$. When does the {\it concatenation} of $f_1(z),e_1(z)$ and $f_2(z),e_2(z)$, i.e., $f,z$ on $[\alpha,\beta]$ defined as
\bq
\begin{array}{l}
f(z)=f_1(z), \quad z \in [\alpha, \gamma), \quad f(z)=f_2(z), \quad z \in (\gamma, \beta] \\[2mm]
e(z)=e_1(z), \quad z \in [\alpha, \gamma), \quad e(z)=e_2(z), \quad z \in (\gamma, \beta]
\end{array}
\eq
satisfy \eqref{InfDir} on the interval $[\alpha, \beta]$ in a {\it weak sense}? This holds if and only if $b_1(\gamma) = b_2(\gamma)$, where
\bq
b_i(z)= Z(\frac{d}{dz}) \ell_i(z), \, \begin{bmatrix} f_i(z) \\[2mm] e_i(z) \end{bmatrix}= \begin{bmatrix} F^\top (\frac{d}{dz}) \\[2mm]  E^\top (\frac{d}{dz})\end{bmatrix} \ell_i (z), i=1,2
\eq
Thus in general $b(\gamma)$ provides exactly the information needed to extend a solution $f(z),e(z)$ on an interval $[\alpha,\gamma]$ to a larger interval $[\alpha,\beta]$. This has been discussed in more detail in the related context of linear partial differential equations involving both spatial and time variables in \cite{helmke}.

\subsection{Differential operator effort constraints}
Another interesting case to be considered concerns differential operator Dirac structures arising from skew-adjoint differential operators and differential operator constraints on the effort variables $e$. This case is well-motivated from an applications point of view. 
Let $J(\frac{d}{dz})$ be a linear differential operator which is formally skew adjoint, i.e., $J(s)=-J^\top(-s)$. Consider the linear space of functions $f(z),e(z)$ satisfying the implicit set of differential equations
\bq
\label{effconstraints}
\begin{array}{l}
\{(f(z),e(z)) \mid  \exists \lambda (z) \mbox{ such that } \\[2mm]
 f(z)= J(\frac{d}{dz})e(z) + G^\top(- \frac{d}{dz})\lambda (z), \; G(\frac{d}{dz}) e(z)=0 \}
\end{array}
\eq
where $G(\frac{d}{dz})$ is a linear differential operator representing the {\it effort constraints} $G(\frac{d}{dz}) e(z)=0$. The vector $\lambda (z)$ represents a vector of Lagrange multiplier functions. For any such $f_i(z),e_i(z), i=1,2,$ belonging to the set \eqref{effconstraints} for some $\lambda_i(z), i=1,2,$ we compute
\bq
\begin{array}{l}
 <(f_1(z),e_1(z)),(f_2(z),e_2(z))> = \\[2mm]
 f_1(z)^\top e_2(z) + e_1(z)^\top f_2(z) = \\[2mm]
 e_2^\top [J(\frac{d}{dz})e_1 + G^\top(- \frac{d}{dz})\lambda_1] + e_1^\top [J(\frac{d}{dz})e_2 + G^\top(- \frac{d}{dz})\lambda_2]  \\[2mm]
 = \,e_2^\top J(\frac{d}{dz})e_1 + e_1^\top J(\frac{d}{dz})e_2 +\\[2mm]
  e_2^\top G^\top (- \frac{d}{dz})\lambda_1 + e_1^\top G^\top (- \frac{d}{dz})\lambda_2
 \end{array}
 \eq
The integral of the term $e_2^\top J(\frac{d}{dz})e_1 + e_1^\top J(\frac{d}{dz})e_2$ over any finite interval $[\alpha,\beta]$ can be computed as follows. Integration by parts and use of $J(s)=-J^\top(-s)$ yields
\bq
\begin{array}{l}
\int_{\alpha}^{\beta} [e_2(z)^\top J(\frac{d}{dz})e_1(z) + e_1(z)^\top J(\frac{d}{dz})e_2(z)] dz = \\[2mm]
[Z_J(\frac{d}{dz})e_1(z)]^\top \Pi_J Z_J(\frac{d}{dz})e_2(z) |^{\beta}_{\alpha},
\end{array}
\eq
where the differential operator $Z_J(\frac{d}{dz})$ and the matrix $\Sigma_J$ are obtained by the two-variable polynomial factorization
\bq
\label{J}
J(\zeta) + J^\top(\eta) = (\zeta + \eta) Z^\top_J(\zeta) \Sigma_J Z(\eta)
\eq
Analogously, the integral of $e_2^\top G^\top (- \frac{d}{dz})\lambda_1 + e_1^\top G^\top (- \frac{d}{dz})\lambda_2$ over $[\alpha, \beta]$ yields by integration by parts
\bq
\begin{array}{l}
\int_{\alpha}^{\beta} [e_2^\top G^\top (- \frac{d}{dz})\lambda_1 + e_1^\top G^\top (- \frac{d}{dz})\lambda_2] dz = \\[2mm]
[Z_G(\frac{d}{dz})e_2(z)]^\top \Pi_G V_G(\frac{d}{dz})\lambda_1(z) + \\[2mm]
[Z_G(\frac{d}{dz})e_1(z)]^\top \Pi_G V_G(\frac{d}{dz})\lambda_2(z) |^{\beta}_{\alpha},
\end{array}
\eq
where the differential operators $Z_G(\frac{d}{dz}), V_G(\frac{d}{dz})$ and the matrix $\Pi_G$ are obtained by the two-variable polynomial factorization
\bq
\label{G}
G^\top(-\eta) - G^\top (\zeta) = (\zeta + \eta) Z^\top_G(\zeta) \Pi_G V_G(\eta)
\eq
Thus two sources for boundary terms are arising. The first one defined by $Z_J(\frac{d}{dz})$ and $\Sigma_J$, resulting from the formally skew-adjoint linear differential operator $J(\frac{d}{dz})$, and the second defined by $Z_G(\frac{d}{dz}), V_G(\frac{d}{dz})$ and $\Pi_G$, resulting from the effort constraints $G(\frac{d}{dz}) e(z)=0$ in \eqref{effconstraints}. Adjoining these boundary terms then yields as before a differential operator Dirac structure on any spatial domain $[\alpha, \beta]$.
\begin{proposition}
Consider a formally skew-adjoint differential operator $J(\frac{d}{dz})$ and a differential operator $G(\frac{d}{dz})$ defining the effort constraints $G(\frac{d}{dz})e(z)=0$. Consider the two variable polynomial factorizations \eqref{J}, \eqref{G}.\\
Then for any $\alpha, \beta \in \mR$ the space defined by \eqref{effconstraints} and the boundary variables
\bq
\begin{array}{l}
b_J(z)= Z_J(\frac{d}{dz}) e(z)\\[2mm]
b_G(z)= Z_G(\frac{d}{dz}) e(z), \, c_G(z)= V_G(\frac{d}{dz}) \lambda (z) 
\end{array}
\eq
defines a Dirac structure with respect to the bilinear form
\bq
\begin{array}{l}
\int_{\alpha}^{\beta}  <(f_1(z),e_1(z)),(f_2(z),e_2(z))>dz  \\[3mm]
\qquad -\,  b^\top_{J1} (\beta) \Sigma_J \, b_{J2}(\beta) \, + \, b^\top_{J1} (\alpha) \Sigma_J \, b_{J2}(\alpha) \\[2mm]
\qquad - \, b^\top_{G2} (\beta) \Pi_G \, c_{G1}(\beta) \, + \, b^\top_{G2} (\alpha) \Pi_G \, c_{G1}(\alpha) \\[2mm]
\qquad - \, b^\top_{G1} (\beta) \Pi_G \, c_{G2}(\beta) \, + \, b^\top_{G1} (\alpha) \Pi_G \, c_{G2}(\alpha)
\end{array}
\eq
\end{proposition}
Furthermore, we conjecture that any differential operator Dirac structure defined by a pair $F(\frac{d}{dz}), E(\frac{d}{dz})$ as above can be also represented as a differential operator Dirac structure defined by a skew-adjoint operator $J(\frac{d}{dz})$ and some $G(\frac{d}{dz})$; analogously to the finite-dimensional case exposed in \cite{DAE1}.

\section{Differential operator Lagrangian subspaces}
Consider again $\mathcal{F} \times \mathcal{E}$, where $\mathcal{E}= \mathcal{F}^*$. In view of applications we will replace the notation $\mathcal{F}=\mR^m$ by $\mathcal{X}=\mR^m$, with $\mathcal{X}$ standing for a linear state space. In fact, in port-Hamiltonian systems theory the flow space $\mathcal{F}$ is actually the {\it tangent space} to the {\it space of energy variables} $\mathcal{X}$. In the present linear case $\mathcal{F}$ obviously can be identified with $\mathcal{X}$. 

Apart from the \emph{symmetric} bilinear form $<\cdot,\cdot>$ as defined in \eqref{plus} there is {\it another} canonically defined bilinear form on $\mathcal{X} \times \mathcal{E}$, defined as
\bq
[(x_1,e_1),(x_2,e_2)]:= <e_1|x_2> - <e_2|x_1>\; ,
\eq
and has the matrix representation
\begin{equation}\label{eq:Je}
J_e = \begin{bmatrix} 0_m & I_m \\[2mm] -I_m & 0_m \end{bmatrix}\; .
\end{equation}
(The skew-symmetric matrix $J_e$ is the matrix representation of the well-known symplectic form.)
As in the case of the bilinear form $<\cdot,\cdot>$, the skew-symmetric bilinear form $[\cdot,\cdot]$ with matrix representation $J_e$ on $\mathcal{X} \times \mathcal{E}= \mathbb{R}^{2m}$ gives rise to the following skew-symmetric bilinear form on the set $\mathcal{C}(\mathbb{R}, \mathcal{X} \times \mathcal{E})$ of smooth functions $(x,e): \mathbb{R} \to \mathcal{X} \times \mathcal{E}$ with compact support:
\begin{equation}\label{-}
\begin{array}{l}
[[ (x_1,e_1), (x_2,e_2)]] := \\[2mm]
 \qquad  \int_{- \infty}^{\infty}  [(x_1(z),e_1(z)),(x_2(z),e_2(z))] dz
 \end{array}
\end{equation}
This is again a skew-symmetric non-degenerate form, in the sense that if $[[ (x_1,e_1), (x_2,e_2)]]=0$ for all compact support $(x_1,e_1)$, then $(x_2,e_2)=0$.
Thus it defines a {\it symplectic form} on $\mathcal{C}(\mathbb{R}, \mathcal{X} \times \mathcal{E})$.

Recall that a subspace $\mathcal{L}$ of a linear space $\mathcal{V}$ with symplectic form $\omega$ is called {\it Lagrangian} if $\mathcal{L}=\mathcal{L}^{\perp}$ where ${}^{\perp}$ denotes the orthogonal complement with respect to the symplectic form $\omega$. Equivalently, $\mathcal{L}$ is Lagrangian if $\omega$ is zero when restricted to $\mathcal{L}$, and moreover $\mathcal{L}$ is maximal with respect to this property. This leads to the following definition of a differential operator Lagrangian subspace over the infinite spatial domain $\mR$. 

Similar to \eqref{InfDir}, consider now systems of differential equations over the infinite spatial domain $\mR$ given by
\bq
\label{InfLag}
P^\top (\frac{d}{dz}) x(z) + S^\top (\frac{d}{dz}) e(z) = 0, \quad z \in \mR,
\eq
where $P(\frac{d}{dz})$ and $S(\frac{d}{dz})$ are square linear differential operators. Denote the corresponding $m \times m$ polynomial matrices by $P(s)$ and $S(s)$ in the indeterminate $s$. Now let $P(s)$ and $S(s)$ satisfy
\bq
\label{sym}
\begin{bmatrix} P^\top (-s) & S^\top (-s) \end{bmatrix} \begin{bmatrix} 0_m & I_m \\[2mm] -I_m & 0_m \end{bmatrix}\begin{bmatrix} P(s) \\[2mm] S(s) \end{bmatrix} =0 , \quad \forall s \in \mC,
\eq
together with the maximal rank assumption
\bq
\label{ranksym}
\rank \begin{bmatrix} P(s) \\[2mm] S(s) \end{bmatrix} = m,  \quad \forall s \in \mC
\eq
It follows again from \cite{rapisarda} that the set $\mathcal{S}_{\mathcal{L}}$ of solutions with compact support of \eqref{InfLag} satisfies
\bq
\int_{- \infty}^{\infty}  [(x_1(z),e_1(z)),(x_2(z),e_2(z))] dz=0
\eq
for all $(x_1,e_1), (x_2,e_2) \in \mathcal{S}_{\mathcal{L}}$. Hence $\mathcal{S}_{\mathcal{L}} \subset \mathcal{S}_{\mathcal{L}}^{\perp}$. By the same reasoning as in the differential operator Dirac structure case it follows that $\mathcal{S}_{\mathcal{L}} = \mathcal{S}_{\mathcal{L}}^{\perp}$, and thus $\mathcal{S}_{\mathcal{L}}$ is an infinite-dimensional Lagrangian subspace.
\begin{proposition}
Consider $P(s),S(s)$ satisfying \eqref{sym} and \eqref{ranksym}. Let $\mathcal{S}_{\mathcal{L}}$ be the set of smooth solutions of \eqref{InfLag} with compact support within the spatial domain $\mR$. Then $\mathcal{S}_{\mathcal{L}}$ is a Lagrangian subspace with respect to the bilinear form \eqref{+}. \end{proposition}

Also, assuming for the moment that the rational matrix $K(s):= - S^{-\top }(s)P^\top (s)$ is {\it proper}, it follows from standard linear system realization theory that the behavior $\mathcal{S}_{\mathcal{L}}$ is generated by a minimal 'state' space system of the form
\bq
\label{state}
\begin{array}{rcl}
\frac{d}{dz} c(z) & = & Ac(z) + Bx(z) \\[2mm]
e(z) & = & Cc(z) + Dx(z) 
\end{array}
\eq
with 'state' vector $c(z) \in \mR^n$, where $n$ is the McMillan degree of $K(s)$. 
Furthermore, by \eqref{sym} the 'transfer' matrix $K(s)$ satisfies 
\bq
K(s)=K^\top (-s)
\eq
Therefore the system \eqref{state} is an input-output Hamiltonian system (\cite{vds14}) in the sense that there exists an invertible matrix $J_i=-J_i^\top $ satisfying
\bq
A^\top J_i + J_i A = 0, \quad B^\top J_i = C, \quad D=D^\top
\eq
In case $K(s)$ is {\it not} proper, the results continue to hold. This follows from the fact that there always exists a partitioning $\{1, \cdots,n\}=I \cup I^c$ such that by defining the polynomial matrix $\widehat{S}(s)$ as the matrix with $k$-th row given by the $k$-th row of $S(s)$ whenever $k \in I$, and equal to the $k$-th row of $P(s)$ whenever $k \in I^c$, and similarly $\widehat{P}(s)$ as the matrix with $k$-th row given by the $k$-th column of $P(s)$ whenever $k \in I$, and equal to the $k$-th row of $S(s)$ whenever $k \in I^c$, then $\widehat{K}(s)= - \widehat{S}^{-\top }(s)\widehat{P}^\top (s)$ is proper. Hence by realization theory we obtain a similar system as in \eqref{state} with inputs $x_k,k\in I, e_k, k \in I^c,$ and outputs $e_k,k\in I, x_k, k \in I^c$.

By two-variable polynomial calculus the skew-symmetric matrix $J_i$ can be directly inferred from the fact that by \eqref{skew} $P^\top (-s)S(s)-S^\top (-s)P(s)=0$, and thus
\bq
\label{facskew}
P^\top (\zeta)S(\eta) - S^\top (\zeta) P(\eta)  = (\zeta + \eta) W^T(\zeta) \Pi W(\eta)
\eq
for some polynomial matrix $W(s)$ and full-rank {\it skew-symmetric} matrix $\Pi$. Hence, {\it without loss of generality} (note that this is fundamentally different from the Dirac structure case!) we can take
\bq
\label{facskew1}
 \Pi= \begin{bmatrix} 0 & I_p \\-I_p & 0 \end{bmatrix} =:J_i
\eq
Then by defining the image representation
\bq
\label{facskew2}
\begin{bmatrix} x \\ e \end{bmatrix} = \begin{bmatrix} P(\frac{d}{dz}) \\ S(\frac{d}{dz}) \end{bmatrix} \ell, \quad \begin{bmatrix} x_{\delta} \\ e_{\delta} \end{bmatrix} = W(\frac{d}{dz}) \ell, 
\eq
we have 
\bq
x_1^\top e_2 - x_2^\top e_1 = \frac{d}{dz} \left(x_{\delta 1}^\top e_{\delta 2} - x_{\delta 2}^\top e_{\delta 1}\right)
\eq
The map $\ell \mapsto Z(\frac{d}{dz}) \ell$ is a {\it minimal} boundary map.

We obtain the following analog of Proposition \ref{prop:dirac}.
\begin{proposition}
Consider $P(\frac{d}{dz}),S(\frac{d}{dz})$ satisfying \eqref{sym}, \eqref{ranksym}. Consider the factorization \eqref{facskew}, \eqref{facskew1}, \eqref{facskew2}.\\
Then for any $\alpha, \beta \in \mR$ the space of $f(z),e(z),b(\alpha),b(\beta)$, which are solutions of 
\bq
\begin{array}{l}
P^\top (\frac{d}{dz}) f(z) + S^\top (\frac{d}{dz}) e(z) = 0, \quad z \in [\alpha, \beta], \\[3mm]
\begin{bmatrix} f(\alpha) \\[2mm] e(\alpha) \end{bmatrix} = \begin{bmatrix} P (\frac{d}{dz}) \\[2mm]  S (\frac{d}{dz})\end{bmatrix} \ell (\alpha), \;
\begin{bmatrix} f(\beta) \\[2mm] e(\beta) \end{bmatrix} = \begin{bmatrix} P (\frac{d}{dz}) \\[2mm]  S (\frac{d}{dz})\end{bmatrix} \ell (\beta)\\[6mm]
\begin{bmatrix} x_{\delta}(\alpha) \\[2mm]  e_{\delta}(\alpha) \end{bmatrix} = W(\frac{d}{dz}) \, \ell (\alpha), \begin{bmatrix} x_{\delta}(\alpha) \\[2mm]  e_{\delta}(\alpha) \end{bmatrix} = W(\frac{d}{dz}) \, \ell (\beta)
\end{array}
\eq
defines a Lagrangian subspace with respect to the bilinear form
\bq
\begin{array}{l}
\int_{\alpha}^{\beta}  [(x_1(z),e_1(z)),(x_2(z),e_2(z))] dz  \\[3mm]
+ \, x_{\delta 1}^\top (\beta) e_{\delta 2}(\beta) -  e_{\delta 1}^\top (\beta) x_{\delta 2}(\beta) \\[3mm]- x_{\delta 1}^\top (\alpha) e_{\delta 2}(\alpha) 
+  \, e_{\delta 1}^\top (\alpha) x_{\delta 2}(\alpha)
\end{array}
\eq
\end{proposition}
A special case of the above is provided by $P(s)=I_m$ and $S(\frac{d}{dz})$ a formally self-adjoint differential operator, i.e.,
\bq
S^\top (-s)=S(s)
\eq
All of this will be further investigated in \cite{journalpaper}. In particular, in \cite{journalpaper} the exposed theory of differential operator Dirac structures will be {\it combined} with the theory of differential operator Lagrangian subspaces in order to give a general definition of {\it port-Hamiltonian systems} defined by differential operators. This combination is based on coupling a differential operator Dirac structure with elements $(f(z,t),e(z,t))$ with a differential operator Lagrangian subspace with elements $(x(z,t),e(z,t))$ by setting $f(z,t)=-\dot{x}(z,t)$; in the spirit of the theory of linear DAE port-Hamiltonian systems exposed in \cite{DAE1}, \cite{beattie}; see also \cite{DAE2} for the nonlinear extension.

\section{Conclusions}
It has been shown how, by using two-variable polynomial matrix calculus, differential operator Dirac structures can be defined on finite intervals, starting from general pairs of linear differential operators defining a Dirac structure on the whole real line. Of particular interest for applications is the closely related class of Dirac structures with boundary variables derived from formally skew-adjoint operators together with differential operator effort constraints. In order that the boundary variables can be split into boundary flow and effort variables an extra condition is identified on the coefficient matrix of a two-variable polynomial matrix. Using a different canonical bilinear form it is shown how the construction can be also used for the definition of Lagrangian subspaces on bounded intervals. In \cite{journalpaper} this will be employed for an extended definition of boundary control port-Hamiltonian systems.

\end{document}